\documentclass[10pt,review]{siamltex}
\usepackage{color}
\usepackage[usenames,dvipsnames,svgnames,table]{xcolor}
\usepackage{hyperref}
\hypersetup{
     colorlinks   = true,
     citecolor    = blue,
     linkcolor    = blue,
     urlcolor     = blue
}
\usepackage{amssymb,amsfonts,amsmath,latexsym,dsfont,bm}
\usepackage{tikz}
\usepackage{pgfplots}
\usetikzlibrary{shapes}
\usepackage{graphicx}
\usepackage{algorithm,algorithmic}
\usepackage{multirow}
\usepackage{bm}
\usepackage{siunitx}

\usepackage{enumitem}

\graphicspath{{./imgArxiv/}}

\newcommand{\rottext}[1]{\rotatebox{90}{\hbox to 30mm{\hss #1\hss}}}
\newcommand{\rottextt}[1]{\rotatebox{90}{\hbox to 15mm{\hss #1\hss}}}
\newcommand{\R}{\mathbb{R}}

\newcommand{\argmin}{\mathrm{argmin}}

\renewcommand{\phi}{\mathbf{\varphi}}

\newcommand{\m}{\textbf{m}}														% transformation	
														 
\newcommand{\bfd}{\mathbf{d}}							 
														
\newcommand{\bP}{\mathbf{P}}

\newcommand{\bfT}{\mathbf{T}}

\newcommand{\bfA}{\mathbf{A}}

\newcommand{\bfP}{\mathbf{P}}
\newcommand{\bfI}{\mathbf{I}}

\newcommand{\bfn}{\mathbf{n}}

\newcommand{\bfx}{\mathbf{x}}

\newcommand{\bfy}{\mathbf{y}}

\newcommand{\bfL}{\mathbf{L}}
\newcommand{\bfw}{\mathbf{w}}

\newcommand{\bfr}{\mathbf{r}}

\newcommand{\bfQ}{\mathbf{Q}}
\newcommand{\bfJ}{\mathbf{J}}
\newcommand{\bfS}{\mathbf{S}}

\newcommand{\x}{x}														% vector in R^d
\renewcommand{\xi}[1]{\x_{#1}}								% i-th component of point x

\newcommand{\hf}{\frac{1}{2}}

\newcommand{\bfK}{\mathbf{K}}

\newcommand*\samethanks[1][\value{footnote}]{\footnotemark[#1]}

\title{LAP: a Linearize and Project Method for Solving Inverse Problems with Coupled Variables}
\author{ James L. Herring\thanks{Department of Mathematics and Computer Science, Emory University, Atlanta, Georgia, USA. ({\tt \{jlherri,nagy,lruthotto\}@emory.edu})} \and James G. Nagy\samethanks[1] \and Lars Ruthotto\samethanks[1]}

\newlength\iwidth
\newlength\iheight

\usepackage{tabularx}
\usepackage{graphicx}

\begin{document}
\maketitle

\begin{abstract}    Many inverse problems involve two or more sets of variables that represent different physical quantities but are tightly coupled with each other. For example, image super-resolution requires joint estimation of the image and motion parameters from noisy measurements. Exploiting this structure is key for efficiently solving these large-scale optimization problems, which are often ill-conditioned. 
    
    In this paper, we present a new method called Linearize And Project (LAP) that offers a flexible framework for solving inverse problems with coupled variables. LAP is most promising for cases when the subproblem corresponding to one of the variables is considerably easier to solve than the other. LAP is based on a Gauss--Newton method, and thus after linearizing the residual, it eliminates one block of variables through projection. Due to the linearization, this block can be chosen freely.
    Further, LAP supports direct, iterative, and hybrid regularization as well as constraints. Therefore LAP is attractive, e.g., for ill-posed imaging problems.
    These traits differentiate LAP from common alternatives for this type of 
problem such as variable projection (VarPro) and block coordinate descent (BCD). 
    Our numerical experiments compare the performance of LAP to BCD and VarPro using three coupled problems whose forward operators are linear with respect to one block and nonlinear for the other set of variables.
\end{abstract}

\begin{keywords}
Nonlinear Least-Squares, Gauss--Newton Method, Inverse Problems, Regularization, Image Processing, Variable Projection
\end{keywords}

\begin{AMS}
    65F10, % numerical analysis -> numerical linear algebra -> iterative methods for linear systems
    65F22, % numerical analysis -> numerical linear algebra -> ill-posednes and regularization
    65M32 % numerical analysis -> PDEs -> inverse problems
\end{AMS}

\pagestyle{myheadings} \thispagestyle{plain} \markboth{J.L.~HERRING, J.G.~NAGY,  AND L.~RUTHOTTO}{LINEARIZE AND PROJECT FOR INVERSE PROBLEMS}

% new section
\section{Introduction} % (fold)
\label{sec:intro}

We present an efficient Gauss--Newton method called Linearize And Project  (LAP) for solving large-scale optimization problems whose variables consist of two or more blocks representing, e.g., different physics.
Problems with these characteristics arise, e.g., when jointly reconstructing image and motion parameters from a series of noisy, indirect, and transformed measurements. LAP is motivated by problems in which the blocks of variables are nontrivially coupled with each other, but some of the blocks lead to well-conditioned and easy-to-solve subproblems. As two examples of such problems arising in imaging, we consider super resolution~\cite{ParkEtAl2003, FarsiuEtAl2004, ChungEtAl2006} and motion corrected Magnetic Resonance Imaging (MRI)~\cite{BatchelorEtAl2005, CorderoEtAl2016}. 

A general approach to solving coupled optimization problems is to use alternating minimization strategies such as Block Coordinate Descent (BCD)~\cite{HardieEtAl1997, NocedalWright1999}. These straightforward approaches can be applied to most objective functions and constraints and also provide flexibility for using various regularization strategies. However, alternating schemes have been shown to converge slowly for problems with tightly coupled blocks~\cite{NocedalWright1999, ChungEtAl2006}.

One specific class of coupled optimization problems, which has received much attention, is separable nonlinear least-squares problems; see, e.g.,~\cite{GolubPereyra1973, GolubPereyra2003, OLearyRust2013}. Here, the variables can be partitioned such that the residual function is linear in one block and nonlinear in the other. For brevity, we will refer to the sets of variables as \emph{linear } and \emph{nonlinear block of variables}, respectively.
 One common method for solving such problems is Variable Projection (VarPro)~\cite{GolubPereyra1973, GolubPereyra2003, OLearyRust2013}. The idea is to derive a nonlinear least-squares problem of reduced size by eliminating the the linear block of variables through projections. VarPro is most effective when the projection, which entails solving a linear least-squares problem, can be computed cheaply and accurately. When the number of linear variables in the problem is large, iterative methods can be used to compute the projection~\cite{FohringEtAl2014}, but iterative methods can become inefficient when the least-squares problem is ill-posed. Further, standard iterative methods do not provide bounds on the optimality, which are needed to ensure well-defined gradients; see our discussion in Sec.~\ref{sub:varpro}. We note that recently proposed methods such as in~\cite{EstrinEtAl2017} can provide that information and are, thus, attractive options in these cases. Another limitation of VarPro is that it is not straightforward to incorporate inequality constraints on the linear variables. Some progress has been made for box-constraints using a pseudo-derivative approach~\cite{SimaVanHuffel2007} leading to approximate gradients for the nonlinear problem. Finally, VarPro limits the options for adaptive regularization parameter selection strategies for the least-squares problem. 

As an alternative to the methods above, we propose the LAP method for solving coupled optimization problems. Our contributions can be summarized as follows:
\begin{itemize}[leftmargin=5mm,itemindent=0mm,itemsep=1mm]
\item We propose an efficient iterative method called LAP that computes the Gauss--Newton step at each iteration by eliminating one block of variables through projection and solving the reduced problem iteratively. Since projection is performed after linearization, any block can be eliminated. Hence the LAP framework offers superior flexibility compared to existing projection-based approaches, e.g., by supporting various types of regularization strategies and the option to impose constraints for all variables. 
\item We demonstrate LAP's flexibility for different regularization strategies including Tikhonov regularization using the discrete gradient operator with a fixed regularization parameter and a hybrid regularization approach~\cite{ChungEtAl2008}, which simultaneously computes the search direction and selects an appropriate regularization parameter at each iteration. 
\item We use projected Gauss--Newton to implement element-wise lower and upper bound constraints with LAP on the optimization variables. This is a distinct advantage over VarPro, where previously proposed methods for incorporating inequality constraints require using approximate gradients.
\item We present numerical experiments for several separable nonlinear least-squares problems. The problems are characterized by linear imaging models and nonlinear motion models with applications including 2D and 3D super-resolution and motion correction for magnetic resonance imaging. We compare the performance of LAP with block coordinate descent (BCD) and variable projection (VarPro) by analyzing convergence, CPU timings, number of matrix-vector multiplications, and the solution images and motion parameters.
\item We provide our MATLAB implementation of LAP including the examples used in this paper freely at
\begin{center}
	\url{https://github.com/herrinj/LAP}
\end{center}

\end{itemize}

Our paper is organized as follows: Sec.~\ref{sec:discrete} introduces a general formulation of the motion-corrected imaging problem, which we use to motivate and illustrate LAP, and briefly reviews BCD and VarPro; Sec.~\ref{sec:LAP} explains our proposed scheme, LAP, for the coupled Gauss--Newton iteration along with a discussion of regularization options and implementation of bound constraints using projected Gauss--Newton; and Sec.~\ref{sec:experiments} provides experimental results for several examples using LAP and compares it with BCD and VarPro. We end with a brief summary and some concluding remarks.

% new section
\section{Motion-Corrected Imaging Problem} % (fold)
\label{sec:discrete}
In this section, we give a general description of coupled optimization problems arising in image reconstruction from motion affected measurements and briefly review Block Coordinate Descent (BCD) and Variable Projection (VarPro). 

We follow the guidelines in~\cite{Modersitzki2009}, and consider images as continuously differentiable and compactly supported functions on a domain of interest, $\Omega \subset \mathbb{R}^d$ (typically, $d = 2$ or $3$). We assume that the image attains values in a field $\mathbb{F}$ where $\mathbb{F}=\mathbb{R}$ corresponds to real-valued and $\mathbb{F}=\mathbb{C}$ to complex-valued images. We denote by $\bfx \in \mathbb{F}^n$ a discrete image obtained by evaluating a continuous image at the cell-centers of a rectangular grid with $n$ cells. 

The discrete transformation $\bfy\in\R^{d\cdot n}$ is obtained by evaluating a function \linebreak$y: \Omega \to \R^d$ at the cell-centers and can be visualized as a transformed grid. For the rigid, affine transformations of primary interest in this paper, transformations are comprised of shifts and rotations and can be defined by a small set of parameters denoted by the variable $\bfw$, but in general, the number of parameters defining a transformation may be large. We observe that under a general transformation $\bfy(\bfw)$, the cell-centers of a transformed grid do not align to the cell-centers of the original grid, so to evaluate a discretized image under a transformation $\bfy$, we must interpolate using the known image coefficients $\bfx$. This interpolation can be represented via a sparse matrix $T(\bfy(\bfw)) \in \mathbb{R}^{n \times n}$ determined by the transformation. For the examples in this paper, we use bilinear or trilinear interpolation corresponding to the dimension of the problem, but other alternatives are possible; see, e.g.,~\cite{Modersitzki2009} for alternatives. The transformed version of the discrete image $\bfx$ is then expressed as a matrix-vector product $T(\bfy(\bfw)) \bfx$. 

Using the above definitions of discrete images and their transformations, we then consider the discrete, forward problem for $N$ distinct data observations, 
\begin{equation}\label{eq:fwd}
    \bfd_k = K_k T(\bfy(\bfw_k))\  \bfx + \epsilon_k, \quad \text{ for all } \quad k=1,2,\ldots,N,
\end{equation}
where $\bfd_k \in \mathbb{F}^{\m_k}$ is the measured data, $K_k \in \mathbb{F}^{m_k \times n}$ is a matrix corresponding to the problem-specific image operator, and $\epsilon_k$ is image noise, which we assume to be independently and identically distributed Gaussian noise. 

In this paper, we focus on the case in which the motion can be modeled by a small number of parameters.
In our case, $\bfw_k \in \R^3$ or $\mathbb{R}^6$ models rigid transformations in 2D and 3D, respectively. The total dimension of the motion parameters across all data measurements is then given by $p = 3N$ or $p=6N$ for $d=2$ and $d=3$, respectively. In the application at hand, we note that $p \ll n$. To simplify our notation, we use the column vectors $\bfd \in \mathbb{F}^{m}$ where $m = m_k \cdot N$ and $\bfw \in \R^p$ to denote the data and motion parameters for all $N$ measurements, respectively.

Given a set of measurements $\{\bfd_1, \bfd_2, \ldots, \bfd_N\}$, the motion-corrected imaging problem consists of jointly estimating the underlying image parameters and the motion parameters in~\eqref{eq:fwd}. We formulate this as the coupled optimization problem
\begin{equation}\label{eq:optProb}
    \begin{split}
        \min_{\bfx \in \mathcal{C}_x, \bfw \in \mathcal{C}_w} & \Phi(\bfx, \bfw) = \hf \| \bfK \bfT(\bfw) \bfx - \bfd \|^2 + \frac{\alpha}{2} \| \bfL \bfx \|_2^2 \text{,} 
    \end{split}
\end{equation}

\noindent where the matrices $\bfK$ and $\bfT$ have the following block structure
\begin{equation*}
\bfK = \begin{bmatrix}
K_1 & & & \\
& K_2 & & \\
& & \ddots & \\
& & & K_N \\
\end{bmatrix} \text{\quad and \quad} 
\bfT(\bfw) = \begin{bmatrix}
T(\bfy(\bfw_1))\\
T(\bfy(\bfw_2))\\
\vdots \\
T(\bfy(\bfw_N))
\end{bmatrix}.
\end{equation*}

\noindent We denote by $\mathcal{C}_x \subset \mathbb{F}^n$ and $\mathcal{C}_w \subset \mathbb{R}^p$ rectangular sets used to impose bound constraints on the image and motion parameters.  Lastly, we regularize the problem by adding a chosen regularization operator, $\bfL$ (e.g., a discrete image gradient or the identity matrix) and a regularization parameter, $\alpha>0$, that balances minimizing the data misfit and the regularity of the reconstructed image. We note that finding a good regularization parameter $\alpha$ is a separate, challenging problem which has been widely researched~\cite{ChungEtAl2008, deSturlerKilmer2011, HaberOldenburg2000, Vogel2002}. One strength of LAP is that, for $\bfL = \bfI$, it allows for regularization methods that automatically select the $\alpha$ parameter~\cite{ChungEtAl2008,GazzolaNagy2014,GazzolaNovati2014,GazzolaEtAl2014}. In our numerical experiments, we investigate one such hybrid method for automatic regularization parameter selection \cite{ChungEtAl2008} as well as direct regularization using a fixed $\alpha$ value.

Problems of the form~\eqref{eq:optProb} are often referred to as separable nonlinear least-squares problems. They are a specific class of coupled optimization problems characterized by being nonlinear in one block of variables, $\bfw$, and linear in the other, $\bfx$. Several optimization approaches exist for solving such problems. One option is a fully coupled Gauss--Newton approach to optimize over both sets of variables simultaneously by solving a single linear system. However, this approach does not exploit the convexity of the problem in the image variables, resulting in small Gauss--Newton steps due to the nonlinearity of the motion parameters \cite{ChungEtAl2006}. Two methods, which have been shown to be preferable to the fully coupled optimization for solving separable nonlinear least-squares problems, are Block Coordinate Descent, which represents a fully decoupled optimization approach, and Variable Projection, which represents a partially coupled, optimization approach. We now provide a brief review of those two methods before introducing LAP.

% new section
\subsection{Block Coordinate Descent (BCD)}

BCD represents a fully decoupled approach to solving coupled optimization problems such as~\eqref{eq:optProb}; see, e.g.,~\cite{NocedalWright1999}. In BCD, the optimization variables are partitioned into a number of blocks. The method then sequentially optimizes over one block of variables while holding all the others fixed. After one cycle in which all subsets of variables have been optimized, one iteration is completed. The process is then iterated until convergence. For this paper, we separate the variables into the two subsets of variables suggested by the structure of the problem, one for the image variables and another for the motion variables. At the $k$th iteration, we fix $\bfw_k$ and obtain the updated image $\bfx_{k+1}$ by solving
\begin{equation}
\label{eq:BCD1}
\bfx_{k+1} = \underset{\bfx \in \mathcal{C}_x}{\argmin} \text{ } \Phi(\bfx, \bfw_k). 
\end{equation}
Afterwards, we fix our new guess for $\bfx_{k+1}$ and optimize over the motion $\bfw$,
\begin{equation}
\label{eq:BCD2}
\bfw_{k+1} = \underset{\bfw \in \mathcal{C}_w}{\argmin} \text{ } \Phi(\bfx_{k+1}, \bfw).
\end{equation}
\noindent These two steps constitute a single iteration of the method. We note that BCD is decoupled in the sense that while optimizing over one set of variables, we neglect optimization over the other. This degrades convergence for tightly coupled problems~\cite{NocedalWright1999}. However, BCD has many advantages. It is applicable to general coupled problems including ones that are nonlinear in all blocks of variables. Also, it allows for straightforward implementation of bound constraints and supports various types of regularization. For our numerical experiments, we solve the BCD imaging problem~\eqref{eq:BCD1} inexactly using a single step of projected Gauss--Newton with bound constraints and various regularizers, which we introduce in Section~\ref{sec:Reg}. The optimization problem in the second step~\eqref{eq:BCD2} is small-dimensional and separable, and we perform a single step of Gauss--Newton with a direct solver to compute the search direction.

% new section
\subsection{Variable Projection (VarPro)}\label{sub:varpro}

VarPro is frequently used to solve separable nonlinear least-squares problems such as the one in~\eqref{eq:optProb}; see, e.g.,~\cite{GolubPereyra2003, OLearyRust2013}. The key idea in VarPro is to eliminate the linear variables (here, the image) by projecting the problem onto a reduced subspace associated with the nonlinear variables (here, the motion) and then solving the resulting reduced, nonlinear optimization problem. In our problem~\eqref{eq:optProb}, eliminating the image variables requires solving a linear least-squares problem involving the matrix $T(\bfw)$ that depends on the current motion parameters. We express the projection by
\begin{equation}\label{eq:VarPro1}
\begin{split}
\bfx(\bfw) &=  \underset{\bfx \in \mathcal{C}_x}{\argmin}\quad  \Phi(\bfx, \bfw).
\end{split}
\end{equation}
\noindent Substituting this expression in for $\bfx$, we then obtain a reduced dimensional problem in terms of the nonlinear variable $\bfw$,
\begin{equation} \label{eq:VarPro2}
\begin{split}
& \min_{\bfw \in \mathcal{C}_w} \Phi(\bfx(\bfw), \bfw).
\end{split}
\end{equation}
\noindent The reduced problem is solved to recover the motion parameters, noting that by solving~\eqref{eq:VarPro1} at each iteration, we simultaneously recover iterates for the image. Assuming $\mathcal{C}_w = \mathbb{R}^p$ (unconstrained case), we see that the first-order necessary optimality condition in~\eqref{eq:VarPro2} is
\begin{equation}\label{eq:VarProOpt}
    0 = \nabla_\bfw \Phi(\bfx(\bfw),\bfw) + \nabla_\bfw \bfx(\bfw) \nabla_\bfx \Phi(\bfx(\bfw), \bfw).
\end{equation}
Note that in the absence of constraints on $\bfx$ (i.e.,  $\mathcal{C}_x = \mathbb{F}^n$) the second term on the right hand side of \eqref{eq:VarProOpt}
vanishes due to the first-order optimality condition of~\eqref{eq:VarPro2}. However, this is not necessarily the case when $\mathcal{C}_x \neq \mathbb{F}^n$ or when~\eqref{eq:VarPro1} is solved with low accuracy. In those cases $\nabla_\bfx \Phi(\bfx(\bfw),\bfw)$ does not equal $0$ and computing $\nabla_\bfw \bfx(\bfw)$, which can be as hard as solving the original optimization problem \eqref{eq:optProb}, is inevitable. In such situations, neglecting the second term in~\eqref{eq:VarProOpt} may considerably degrade the performance of VarPro.

% new section
\section{Linearize and Project (LAP)}
\label{sec:LAP}

We now introduce the LAP method for solving the coupled optimization problem~\eqref{eq:optProb}. We begin by linearizing the residual in~\eqref{eq:optProb} following a standard Gauss--Newton framework. In each iteration, computing the search direction then requires solving a linear system that couples the image and motion parameters. We do this by projecting the coupled problem onto one block of variables. This offers flexibility in terms of regularization and can handle bound constraints on both the motion and image variables. Furthermore, it also allows the user to freely choose which block of variables is eliminated via projection. 

We introduce our approach to solving the linear coupled problem for the Gauss--Newton step by breaking the discussion into several subsections. We start with a subsection introducing our strategy of projection onto the image space for an unconstrained problem where $\mathcal{C}_x = \mathbb{F}^n$ and $\mathcal{C}_w = \mathbb{R}^p$. This is followed by a subsection on the various options for image regularization that our projection approach offers. Lastly, we extend LAP to a projected Gauss--Newton framework to allow for element-wise bound constraints on the solution, i.e., when $\mathcal{C}_x$ and $\mathcal{C}_w$ are proper subsets of $\mathbb{F}^n$ and $\mathbb{R}^p$, respectively.

% new section
\subsection{Linearizing the Problem}
\label{sec:Lin}

We begin considering a Gauss--Newton framework to solve the coupled problem~\eqref{eq:optProb} in the unconstrained case, i.e., $\mathcal{C}_x = \mathbb{F}^n$ and $\mathcal{C}_w = \mathbb{R}^p$. To solve for the Gauss--Newton step at each iteration, we first reformulate the problem by linearizing the residual $\bfr(\bfx,\bfw) = \bfK \bfT(\bfw) \bfx - \bfd$ around the current iterate, $(\bfx_0,\bfw_0)$. Denoting this residual as $\bfr_0 = \bfr(\bfx_0, \bfw_0)$, we can write its first-order Taylor approximation as
\begin{equation}\label{eq:linRes}
\bfr (\bfx_0 + \delta\bfx, \bfw_0 + \delta \bfw) \approx \bfr_0 + \begin{bmatrix}
\bfJ_x & \bfJ_w
\end{bmatrix} \begin{bmatrix}
\delta \bfx \\ \delta \bfw  
\end{bmatrix},
\end{equation}
\noindent where $\bfJ_x = \nabla_x \bfr_0^\top$ and $ \bfJ_w = \nabla_w \bfr_0^\top$  are the Jacobian operators with respect to the image and motion parameters, respectively. They can be expressed as
\begin{align*}
     \bfJ_x  &=  \bfT^\top \bfK^\top, \quad \text{ and }\quad \bfJ_w = {\rm diag}(\nabla_{\bfw_1}\left(T(\bfy(\bfw_1))\bfx\right), \ldots,\nabla_{\bfw_N}\left(T(\bfy(\bfw_N))\bfx\right) )^\top \bfK^\top,
 \end{align*}
\noindent where each  term $\nabla_{\bfw_k}\left( T(\bfy(\bfw_k)) \bfx \right)$ is the gradient of the transformed image at the transformation $\bfy(\bfw_k)$; see~\cite{ChungEtAl2006} for a detailed derivation. Both $\bfJ_x \in \mathbb{F}^{m \times n} $ and $\bfJ_w \in \mathbb{F}^{m \times p}$ are sparse and we recall that in the application at hand $p \ll n$. 

After linearizing the residual around the current iterate, we substitute the approximation~\eqref{eq:linRes} for the residual term in~\eqref{eq:optProb} to get  
\begin{equation}\label{eq:optProb2}
\begin{split}
         \min_{\delta \bfx, \delta \bfw} &   \hat{\Phi}(\delta \bfx,\delta \bfw) =  \frac{1}{2} \left\lVert \bfJ_x \delta \bfx + \bfJ_w \delta \bfw + \bfr_0 \right\rVert^2 + \frac{\alpha}{2} \|\bfL(\bfx_0 + \delta \bfx) \|^2  \\
\end{split}
\end{equation}
\noindent By solving this problem we obtain the updates for the image and motion parameters, denoted by $\delta \bfx$ and  $\delta \bfw$, respectively. As~\eqref{eq:optProb2} is based on a linearization it is common practice to solve it only to a low accuracy; see, e.g.,~\cite{NocedalWright1999}. Note that solving \eqref{eq:optProb2} directly using an iterative method equates to the fully coupled Gauss--Newton approach mentioned in Sec.~\ref{sec:discrete}, which has been observed to converge slowly during optimization. This motivates LAP's projection-based strategy to solving the linearized problem for the Gauss--Newton step, which we now present.

% new section
\subsection{Projecting the problem onto the image space}
\label{sec:Proj}

Recall that VarPro is restricted to projecting onto the nonlinear variables (in our case, this would be the motion parameters, which requires us to solve one large-scale image reconstruction problem per iteration). With our framework, there is no such restriction; the coupled linear problem \eqref{eq:optProb2} can be projected onto either set of variables. Because there are a small number of nonlinear variables, to solve the coupled linear problem in~\eqref{eq:optProb2}, we propose projecting the problem onto the image space and solving for $\delta \bfx$. To project, we first note that the first-order optimality condition for the linearized problem with respect to $\delta \bfw$ is 
$$
0 = \nabla_{\delta \bfw} \big( \bfJ_x \delta \bfx + \bfJ_w \delta \bfw + \bfr_0 \big)^\top \big( \bfJ_x \delta \bfx + \bfJ_w \delta \bfw + \bfr_0 \big) 
$$
or equivalently,
$$
0 = 2 \big( \bfJ_w^{\top} \bfJ_w \delta \bfw + \bfJ_w^{\top} \bfJ_x \delta \bfx + \bfJ_w^{\top}  \bfr_0 \big).
$$
\noindent Solving this condition for $\delta \bfw$, we get
\begin{equation} \label{eq:deltaw}
\delta \bfw = -\big(\bfJ_w^\top \bfJ_w \big)^{-1} \big( \bfJ_w^\top \bfJ_x \delta \bfx +  \bfJ_w^\top \bfr_0 \big).
\end{equation}
\noindent We can then substitute this expression for $\delta \bfw$ into~\eqref{eq:optProb2} and group terms. This projects the problem onto the image space and gives a new problem in terms of $\delta \bfx$, the Jacobians $\bfJ_x$ and $\bfJ_w$, and the residual $\bfr_0$ given by
$$
\min_{\delta \bfx} \hf \left\lVert \left( \bfI - \bfJ_w \big(\bfJ_w^\top \bfJ_w \big)^{-1} \bfJ_w^\top \right) \bfJ_x \delta \bfx + \left(\bfI - \bfJ_w \big(\bfJ_w^\top \bfJ_w \big)^{-1} \bfJ_w^\top \right)\bfr_0 \right\rVert^2 + \frac{\alpha}{2} \|\bfL (\bfx_0 + \delta \bfx) \|^2
$$
or more succinctly,
\begin{equation} \label{eq:probX}
\min_{\delta \bfx} \hf \left\lVert \bfP_{\bfJ_w}^\perp \big( \bfJ_x \delta \bfx + \bfr_0\big) \right\rVert^2 + \frac{\alpha}{2} \|\bfL (\bfx_0 + \delta \bfx) \|^2 
\end{equation}

\noindent where $\bfP_{\bfJ_w}^\perp = \bfI - \bfJ_w (\bfJ_w^\top \bfJ_w)^{-1} \bfJ_w^\top$ is a projection onto the orthogonal complement of the column space of $\bfJ_w$. This least-squares problem can be solved using an iterative method with an appropriate right preconditioner~\cite{Bjorck1996, HestenesStiefel1952, Saad2003}. In particular, we observe that
$$ \bfP_{\bfJ_w}^\perp \bfJ_x = \bfJ_x - \bfJ_w \big(\bfJ_w^\top \bfJ_w \big)^{-1} \bfJ_w^\top \bfJ_x $$
\noindent is  a low rank perturbation of the operator $\bfJ_x$ since ${\rm rank}(\bfP_{\bfJ_w}^\perp) \leq p \ll n$.
Hence, a good problem-specific preconditioner for $\bfJ_x$ should be a suitable preconditioner for the projected operator. 

It is important to emphasize that in our approach, the matrix $\bfJ_w^\top \bfJ_w  \in \mathbb{R}^{p \times p}$ is moderately sized and symmetric positive-definite if $\bfJ_w$ is full rank. Therefore, it is computationally efficient to compute its Cholesky factors once per outer Gauss--Newton iteration and reuse them to invert the matrix when needed. Furthermore, one can use a thin QR factorization of $\bfJ_w$ to compute the Cholesky factors and avoid forming $\bfJ_w^\top \bfJ_w$ explicitly. We use this strategy to increase the efficiency of iteratively solving~\eqref{eq:probX} when using matrix-free implementations, and in practice have not seen rank-deficiency. After solving the preconditioned projected problem for $\delta \bfx$, one obtains the accompanying motion step $\delta \bfw$ via~\eqref{eq:deltaw}.

% new section
\subsection{Regularization} \label{sec:Reg}
After linearizing and projecting onto the space of the linear, image variables, we solve the regularized least-squares problem~\eqref{eq:probX}. For our applications, this problem is high-dimensional, and we use an iterative method to approximately solve it with low accuracy. Also, this problem is ill-posed for our numerical tests and thus requires regularization. Here we discuss the types of direct, iterative, and hybrid methods for regularization that can be used in LAP.
 
We first note that~\eqref{eq:probX} is a Tikhonov least-squares problem. Problems of this form are well-studied in the literature; see, e.g.,~\cite{Hansen1998, EnglEtAl2000, Vogel2002, Hansen2010}. The quality of an iterative solution for this type of problem depends on the selection of an appropriate regularization parameter $\alpha$ and a regularization operator $\bfL$. Common choices for $\bfL$ include the discretized gradient operator $\bfL = \nabla_h$ using forward differences and the identity $\bfL = \bfI$. Additionally,  problems including non-quadratic regularizers, e.g., total variation~\cite{RUDIN:1992kn} or $p$-norm based regularizers can be addressed by solving a sequence of least-squares problems that include weighted $\ell_2$ regularizers \cite{RodrWohl2006,RodrWohl2008,RodrWohl2007}, or more efficiently via a hybrid Krylov subspace approach \cite{GazzolaNagy2014}. Efficient solvers for optimization problems with quadratic regularizers are also a key ingredient of splitting-based methods, e.g., the Split Bregman method for $\ell_1$ regularized problems~\cite{GoldsteinOsher2009}. Thus, while we restrict ourselves to quadratic regularization terms for the scope of this paper, LAP is suitable for a broader class or regularization options. For more information on regularization parameter selection, see \cite{Hansen1998, Vogel2002}.

Another approach is iterative regularization \cite{EnglEtAl2000, Hansen1998}, which aims to stop an iterative method at the iteration that minimizes the reconstruction error. It exploits the fact that the early iterations of iterative methods like Landweber and Krylov subspace methods contain the most important information about the solution and are less affected by noise in the data. In contrast, later iterations of these methods contain less information about the solution and are more noise-affected. This results in a reduction of the reconstruction error for the computed solution during the early iterations of the iterative method followed by its increase in later iterations, a phenomenon known as semi-convergence. However, in practice iterative regularization can be difficult as the reconstruction error at each iterate is unknown and the quality of the solution is sensitive to an accurate choice of stopping iterate.

Hybrid regularization methods represent further alternatives that seek to combine the advantages of direct and iterative regularization, see e.g., \cite{ChungEtAl2008,GazzolaNagy2014,GazzolaNovati2014,GazzolaEtAl2014} and the references therein. Hybrid methods use Tikhonov regularization with a new regularization parameter $\alpha_k$ at each step of an iterative Krylov subspace method such as LSQR~\cite{PaigeSaunders1982}. Direct regularization for the early iterates of a Krylov subspace method is possible due to the small dimensionality of the Krylov subspace. This makes SVD-based parameter selection methods for choosing $\alpha_k$ computationally feasible at each iteration. The variability of $\alpha_k$ at each iteration helps to stabilize the semi-convergence behavior of the iterative regularization from the Krylov subspace method, making it less sensitive to the choice of stopping iteration. Thus, hybrid methods combine the advantages of both direct and iterative regularization while avoiding the cost of SVD-based parameter selection for the full-dimensional problem and alleviating the sensitivity to stopping criteria which complicates iterative regularization.

Each of the methods in this paper necessitates solving a regularized system in the image variable given by~\eqref{eq:probX} for LAP,~\eqref{eq:VarPro1} for VarPro, and~\eqref{eq:BCD1} for BCD. For these problems, we use both direct and hybrid regularization when possible to demonstrate the flexibility of the LAP approach. We begin by running LAP with the discrete gradient regularizer, $\bfL = \nabla_h$ and a fixed $\alpha$. This approach is also feasible within the VarPro framework and is straightforward to implement using BCD. To test hybrid regularization, we use the \texttt{HyBR} method~\cite{ChungEtAl2008} (via the interface given in the IR Tools package~\cite{IRToolsv1}) for LAP and BCD to automatically choose $\alpha$ using a weighted GCV method. We note that because such hybrid regularization methods are tailored to linear inverse problems, they cannot be used to directly solve the nonlinear VarPro optimization problem \eqref{eq:VarPro1}.

% new section
\subsection{Optimization} \label{sec:Optim}
In Sec.~\ref{sec:Lin} and~\ref{sec:Proj}, we introduced LAP for \eqref{eq:optProb} for the case of an unconstrained problem. In practice for imaging problems such as~\eqref{eq:optProb}, bounds on the image and/or motion variables are often known, in which case imposing such prior knowledge into the inverse problem is desirable. This section details how the LAP strategy can be coupled with projected Gauss--Newton to impose element-wise bound constraints on the solutions for $\bfx$ and $\bfw$. We introduce projected Gauss--Newton following the description in~\cite{Haber2014}. The method represents a compromise between a full Gauss--Newton, which converges quickly when applied to the full problem and projected gradient descent, which allows for the straightforward implementation of bound constraints. For projected Gauss--Newton, we separate the step updates $\delta \bfx$  and $\delta \bfw$ into two sets: the set of variables for which the bound constraints are inactive (the inactive set) and the set of variables for which the bound constraints are active (the active set). We denote these subsets for the image by $\delta \bfx_{\mathcal{I}} \subset \delta \bfx$ and $\delta \bfx_{\mathcal{A}} \subset \delta \bfx$ where the subscript $\mathcal{I}$ and $\mathcal{A}$ denote the inactive and active sets, respectively. Identical notation is used for the inactive and active sets for the motion.

On the inactive set, we take the standard Gauss--Newton step at each iteration using the LAP strategy described in Sec.~\ref{sec:Proj}. This step is computed by solving \eqref{eq:optProb2} restricted to the inactive set. Thus, \eqref{eq:probX} becomes
\begin{equation} \label{eq:inactProbX}
\begin{split}
 &\min_{\delta \bfx} \hf \left\lVert \hat{\bfP}_{\bfJ_w}^\perp \big( \hat{\bfJ}_x \delta \bfx_{\mathcal{I}} + \bfr_0\big) \right\rVert^2 + \frac{\alpha}{2} \|\hat{\bfL} (\bfx_{0,\mathcal{I}} + \delta \bfx_{\mathcal{I}}) \|^2, 
\end{split}
\end{equation}
\noindent where $\hat{\bfJ}_x$, $\hat{\bfJ}_w$, $\hat{\bfP}_{\bfJ_w}^\perp$, and $\hat{\bfL}$ are $\bfJ_x$, $\bfJ_w$, $\bfP_{\bfJ_w}$, and $\bfL$ restricted to the inactive set via projection. We then obtain the corresponding motion step on the inactive step by
\begin{equation} \label{eq:inactDeltaW}
\delta \bfw_{\mathcal{I}} = -\big(\hat{\bfJ}_w^\top \hat{\bfJ}_w \big)^{-1} \big( \hat{\bfJ}_w^\top \hat{\bfJ}_x \delta \bfx_{\mathcal{I}} +  \hat{\bfJ}_w^\top \bfr_0 \big).
\end{equation}
\noindent This equation is analogous to \eqref{eq:deltaw} for the unconstrained problem. Note that LAP's projection strategy is only used on the inactive set. Thus, the constraints do not affect the optimality condition for the projection to eliminate the $\delta \bfw_{\mathcal{I}}$ block of variables. The projected least-squares problem~\eqref{eq:inactProbX} is also unaffected. Also, for the special case when the upper and lower bounds on all variables are $-\infty$ and $\infty$ and all variables belong to the inactive set at each iteration, the method reduces to standard Gauss--Newton as presented in Sec.~\ref{sec:Proj}. 

For the active set, we perform a scaled, projected gradient descent step given by
\begin{equation}\label{eq:projGradDes}
\begin{bmatrix}
\delta \bfx_{\mathcal{A}} \\
\delta \bfw_{\mathcal{A}}
\end{bmatrix} =
-\begin{bmatrix}
\tilde{\bfJ}_x^\top \bfr_0 \\ 
\tilde{\bfJ}_w^\top \bfr_0
\end{bmatrix}  
- \alpha \begin{bmatrix}
\tilde{\bfL}^\top \tilde{\bfL} (\bfx_{0,\mathcal{A}} + \delta \bfx_{\mathcal{A}}) \\
{\bf 0}
\end{bmatrix}.
\end{equation}
\noindent where again $\tilde{\bfJ}_x$, $\tilde{\bfJ}_w$, and $\tilde{\bfL}$ represent the projection of $\bfJ_x$, $\bfJ_w$, and $\bfL$ onto the active set. We note that the regularization parameter $\alpha$ should be consistent for both \eqref{eq:inactProbX} and \eqref{eq:projGradDes}. When using direct regularization with a fixed $\alpha$, this is obvious. However, for the hybrid regularization approach discussed in Sec.~\ref{sec:Reg}, a choice is required. We set $\alpha$ on the active set at each iteration to be the same as the $\alpha$ adaptively chosen by the hybrid regularization on the inactive set at the same iterate.

The full step for a projected Gauss--Newton iteration is then given by a scaled combination of steps on the inactive and active sets
\begin{equation}\label{eq:projGN}
\begin{bmatrix}
\delta \bfx \\
\delta \bfw
\end{bmatrix}
 = \begin{bmatrix}
\delta \bfx_{\mathcal{I}} \\
\delta \bfw_{\mathcal{I}}
\end{bmatrix}
 + \gamma \begin{bmatrix}
 \delta \bfx_{\mathcal{A}} \\
 \delta \bfw_{\mathcal{A}}
 \end{bmatrix}.
\end{equation}
\noindent Here, the parameter $\gamma > 0$ is a weighting parameter to reconcile the difference in scales between the Gauss--Newton and gradient descent steps. To select this parameter, we follow the recommendation of ~\cite{Haber2014} and use
\begin{equation}
\label{eq:gamma}
\gamma = \frac{\max \left( \| \delta \bfx_{\mathcal{I}} \|_{\infty}, \| \delta \bfw_{\mathcal{I}} \|_{\infty} \right)}{\max \left( \| \delta \bfx_{\mathcal{A}} \|_{\infty}, \| \delta \bfw_{\mathcal{A}} \|_{\infty} \right)}.
\end{equation}
\noindent This choice of $\gamma$ ensures that the projected gradient descent step taken on the active set is no larger than the Gauss--Newton step taken on the inactive set, and we have used it with no ill effects in practice. 

After combining the steps for both the inactive and active sets using~\eqref{eq:projGN}, we use a projected Armijo line search for the combined step to ensure the next iterate obeys the problem's constraints. A standard Armijo line search chooses a step size $0 < \eta \leq 1$ by backtracking from the full Gauss--Newton step ($\eta = 1$) to ensure a reduction of the objective function \cite{NocedalWright1999}. The projected Armijo line search satisfies a modified Armijo condition given by

\begin{equation}
\label{eq:projArmijo}
\Phi \big( \bP_{\mathcal{C}_\bfx}(\bfx + \eta \delta \bfx),  \bP_{\mathcal{C}_\bfw}(\bfw + \eta \delta \bfw)\big)  \leq \Phi(\bfx, \bfw) + c \eta \bfQ \big( \nabla \Phi(\bfx,\bfw)\big)^\top  \begin{bmatrix} \delta \bfx \\ \delta \bfw \end{bmatrix} . 
\end{equation}

\noindent Here, $\bP_{\mathcal{C}_\bfx}$ and $\bP_{\mathcal{C}_\bfw}$ are projections onto the feasible set for the image and motion variables, respectively, and  $\bfQ \big( \nabla \Phi(\bfx, \bfw)\big)$ is the projected gradient. The constant $c \in (0,1)$ determines the necessary reduction for the line search; we set $c = \num{e-4}$ as suggested in \cite{NocedalWright1999}. Under the projections $\bP_{\mathcal{C}_\bfx}$ and $\bP_{\mathcal{C}_\bfw}$, variables that would leave the feasible region are projected onto the boundary and join the active set for the next iteration. However, the projection does not prevent variables from leaving the active set and joining the inactive set. This necessitates updating the inactive and active sets after the line search at each next iteration. Also note that for the special, unconstrained case, the line search reverts to the standard Armijo line search with no need for projection. However, for this case, the problem reverts to the standard Gauss--Newton framework which should remove the necessity for a line search in most cases.

Lastly, we discuss the choice of stopping criteria for the projected Gauss--Newton method, which again depends on the type of regularization used. For a fixed regularization parameter $\alpha$, we monitor the relative change of the objective function and the norm of the projected gradient including the regularizer term. The projected gradient is the first-order optimality condition of the constrained problem in~\eqref{eq:optProb} and can be computed via a projection~\cite{Beck2014}. The projection is necessary because the bound constraints prevent gradient entries corresponding to variables with minima outside the feasible region from converging to zero. When using hybrid regularization, we must consider the variability of $\alpha_k$ at each Gauss--Newton iteration. Selecting a different $\alpha_k$ parameter at each iteration changes the weight of the regularization term in the objective function~\eqref{eq:optProb} and its projected gradient at each iteration. This makes the full objective function value and projected gradient for these methods unreliable stopping criteria. Instead, we monitor both the norm of the difference between the current iterate and the previous iterate and the difference between the previous and current objective function values for the data misfit term of the objective function. We stop the Gauss--Newton method when either of those values drops below a certain threshold, indicating a stagnation of the method. This allows us to monitor the behavior of the Gauss--Newton iteration without being subject to the fluctuations associated with the varying $\alpha_k$.

% new section
\subsection{Summary of the Method}
We summarize the discussion of LAP by presenting the entire projected Gauss--Newton algorithm using the LAP framework. This provides a framework to efficiently solve the coupled imaging problems of interest while including the flexible options for regularization and bound constraints which helped motivate our approach. The complete algorithm is given in Algorithm~\ref{LAPalgorithm}.
\begin{algorithm}[t]
\caption{Linearize and Project (LAP)} \label{LAPalgorithm}
\begin{algorithmic}[1]
\STATE Given $\bfx_0$ and $\bfw_0$ 
\STATE Compute active and inactive sets, $\mathcal{A}$ and $\mathcal{I}$
\STATE Evaluate $\Phi(\bfx_0,\bfw_0)$, $\bfr_0$, $\bfJ_x^{(0)}$ and $\bfJ_w^{(0)}$
\FOR{$k = 1,2,3,\ldots$}
	\STATE Compute the step the inactive set with LAP using \eqref{eq:inactProbX} and \eqref{eq:inactDeltaW} 
    	\STATE Compute the step on the active set using projected gradient descent \eqref{eq:projGradDes}
   	 \STATE Combine the steps using \eqref{eq:projGN} and \eqref{eq:gamma}
      	 \STATE Perform the projected Armijo line search satisfying \eqref{eq:projArmijo}, update $\bfx_k, \bfw_k$
       	 \STATE Update active and inactive sets, $\mathcal{A}$ and $\mathcal{I}$
    	 \STATE Evaluate $\Phi(\bfx_k, \bfw_k)$,  $\bfr_k$, and $\bfJ_x^{(k)}$ and $\bfJ_w^{(k)}$
    	 \IF{ Stopping criteria satisfied}
        		\RETURN{$\bfx_k, \bfw_k$} 
    	\ENDIF    
\ENDFOR
\end{algorithmic}
\end{algorithm}
%

% new section
\section{Numerical Experiments} % (fold)
\label{sec:experiments}

We now test LAP for three coupled imaging problems. These are a two-dimensional super-resolution problem, a three-dimensional super-resolution problem, and a linear MRI motion correction problem. For all examples, we compare the results using LAP to those of VarPro and BCD. We look at the quality of the resultant image, the ability of the methods to correctly determine the motion parameters across all data frames, the number of iterations required during optimization, the cost of optimization in terms matrix-vector multiplications, and the CPU time to reach a solution. To compare the quality of the resultant image and motion, we use relative errors with respect to the known, true solutions. The matrix-vector multiplications of interest are those by the Jacobian operator associated with the linear, imaging variable $\bfJ_x$. Matrix-vector multiplications by this operator are required in the linear least-squares systems for all three methods and are the most computationally expensive operation within the optimization.

% new section
\subsection{Two-Dimensional Super Resolution} % (fold)
\label{sub:two_dimensional_example}

Next, we run numerical experiments using a relatively small two-dimensional super resolution problem. To construct a super resolution problem with known ground truth image and motion parameters, we use the 2D MRI dataset provided in FAIR~\cite{Modersitzki2009} (original resolution $128\times 128$) to generate 32 frames of low-resolution test data (resolution $32\times 32$) after applying 2D rigid body transformations with randomly chosen parameters. Gaussian white noise is added using the formula
\begin{equation*}
\begin{array}{lll}
\bfd_k = \bar{\bfd}_k + \epsilon_k & \text{ where } & \epsilon_k = \mu \frac{\|\bar{\bfd}_k \|_2}{\| \bfn_k \|_2} \bfn_k
\end{array} 
\end{equation*}
\noindent where $\bar{\bfd_k}$ denotes a noise free low-resolution data frame, $\mu$ is the percentage of noise and $\bfn_k$ is a vector of normally distributed random values with mean $0$ and standard deviation $1$. Our experiments show results with $\mu=1$\%, $2$\%, and $3$\% noise added to the low-resolution data frames. The resulting super-resolution problem is an optimization problem of the form \eqref{eq:optProb}. Here, the imaging operator $\bfK$ is a block diagonal matrix composed of $32$ identical down-sampling matrices $K = K_k$ along the diagonal, which relate the high-resolution image to the lower resolution one via block averaging. The total number of parameters in the optimization is $16,528$ corresponding to $16,384$ for the image and $96$ for the motion.

As mentioned in~\cite{ChungEtAl2006}, the choice of initial guess is crucial in super-resolution problems, so we perform rigid registration of the low-resolution data onto the first volume (resulting in 31 rigid registration problems) to obtain a starting guess for the motion parameters. The resulting relative error of the parameters is around $2$\%. Using these parameters, we solve the linear image reconstruction problem to obtain a starting guess for the image with around $5$\% relative error. 

We then solve the super-resolution problem using LAP, VarPro, and BCD. For all three approaches, we compare two regularization strategies. All three methods are run using the gradient operator, $\bfL = \nabla_h$ with a fixed regularization parameter $\alpha = 0.01$. LAP and BCD are then run with the Golub-Kahan hybrid regularization approach detailed in Section~\ref{sec:Reg} (denoted as \texttt{HyBR} in tables and figures). As previously mentioned, \texttt{HyBR} cannot be applied directly to the VarPro optimization problem \eqref{eq:VarPro1}, so for VarPro we use a fixed $\alpha = 0.01$ and the identity as our second regularizer, $\bfL = \bfI$. To allow for comparison, we use the same regularization parameter as in \cite{ChungEtAl2006}, which may not be optimal for any of the methods presented. For more rigorous selection criteria, we refer the reader to the methods mentioned in Sec.~\ref{sec:Reg}. For LAP and BCD, we add element-wise bound constraints on the image space in the range $[0,1]$  for both choices of regularizer, with both bounds active in practice. The number of active bounds varies for different noise levels and realizations of the problem but can include as many as $30$ -- $35$\% of the image variables for both LAP and BCD. VarPro is run without bound constraints. Neither constraints nor regularization are imposed on the motion parameters for the three methods. 

All three methods require solving two different types of linear systems, one associated with the image variables and another with the motion variables. LAP solves these systems to determine the Gauss--Newton step. We use LSQR~\cite{PaigeSaunders1982} with a stopping tolerance of \textcolor{blue}{$\num{e-2}$} to solve~\eqref{eq:probX} for both regularization approaches. The motion step~\eqref{eq:deltaw} is computed using the Cholesky factors of $\bfJ_w^\top \bfJ_w$. VarPro requires solving the linear system~\eqref{eq:VarPro1} within each function call. For both choices of regularization, we solved this system by running LSQR for a fixed number of $20$ iterations. Recall that this system must be solved to a higher accuracy than the similarly sized systems in LAP and BCD to maintain the required accuracy in the gradient; see Sec.~\ref{sub:varpro}. Gauss--Newton on the reduced dimensional VarPro function~\eqref{eq:VarPro2} then requires solving the reduced linear system in the motion parameters, which we solve using Cholesky factorization on the normal equations. For BCD,  coordinate descent requires alternating solutions for the linear system in the image,~\eqref{eq:BCD1}, and a nonlinear system in the motion parameters,~\eqref{eq:BCD2}. For both of these, we take a single projected Gauss--Newton step. For the image, this is solved using LSQR with a stopping tolerance of \textcolor{blue}{$\num{e-2}$} using \texttt{MATLAB}'s \texttt{lsqr} function for direct regularization and \texttt{HyBR}'s LSQR for hybrid regularization. LSQR with this tolerance is also used for the corresponding problems in the 3D super-resolution and MRI motion correction examples in Sec.~\ref{sub:three_dimensional_example} and~\ref{sub:MRI_motion_example}. For the motion, we use Cholesky factorization on the normal equations. During these solves, we track computational costs for the number of matrix-vector multiplications by the Jacobian operator associated with the image, $\bfJ_x$. For LAP and BCD, these multiplications are required when solving the system for the Gauss--Newton step for the image step, while for VarPro, they are necessary for the least-squares solve within the objective function.

\begin{figure}[t]
    \begin{center}
        \begin{tabular}{c}
        \includegraphics[width=.95\textwidth]{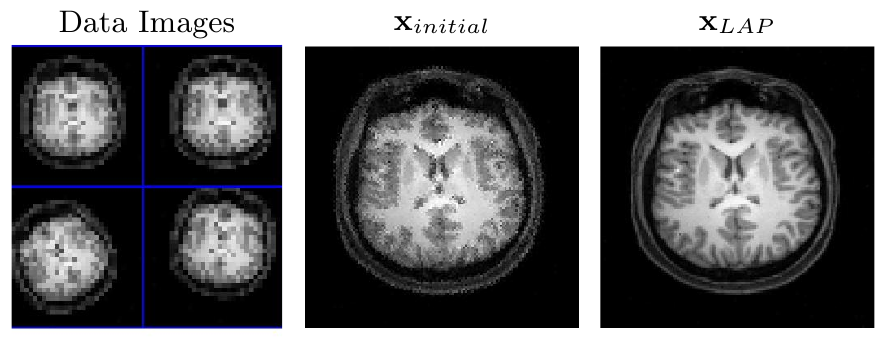} 
        \end{tabular}
    \end{center}
    \caption{The images from left to right show a montage of data frames, the initial guess $\bfx_{initial}$, and the reconstructed image $\bfx_{LAP}$ for LAP using the discrete gradient regularizer for the 2D super resolution problem with 2\% noise.}
    \label{fig:2D_SuperRes_Images}
\end{figure}

\begin{figure}[t]
    \begin{center}
        \includegraphics[width=\textwidth]{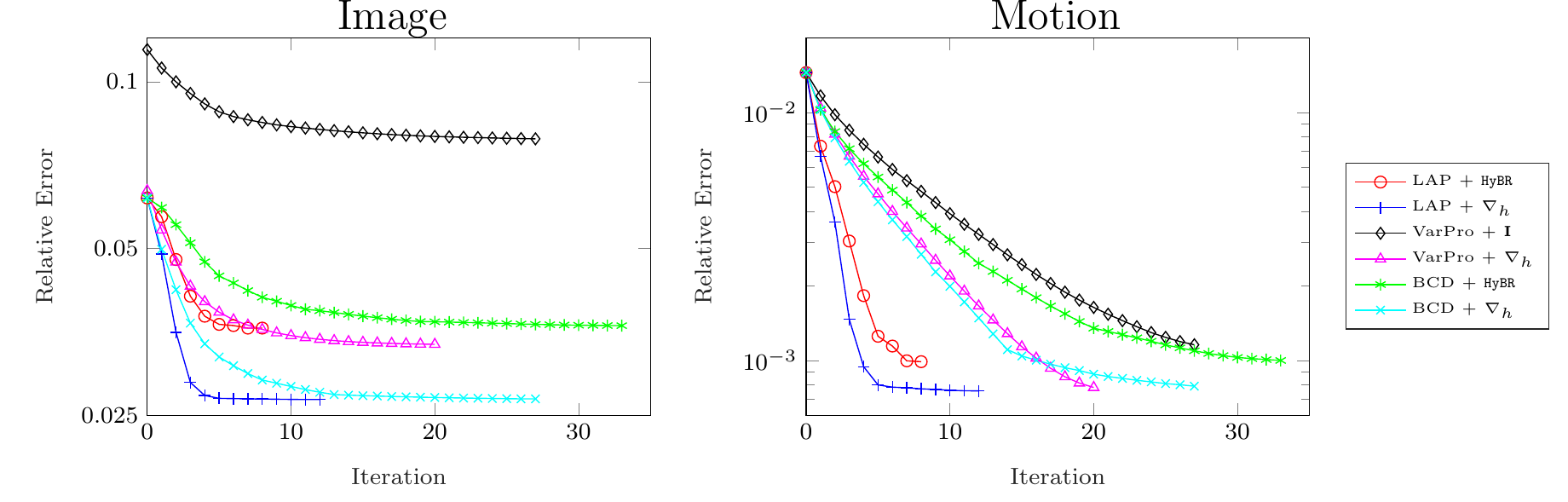}
    \end{center}
    \caption{These plots show the relative errors for both the reconstructed image and the motion parameters for the 2D super resolution problem with $2$\% added noise. We note that methods using the discrete gradient regularizer reach lower minima for the image in this problem, and LAP with the discrete gradient regularizer outperforms both VarPro and BCD in recovering the motion parameters in the early iterations of the method for all noise levels tested.}
    \label{fig:2D_SuperRes_RE}
\end{figure}

We compare results for the methods using the relative errors of the resultant image and motion parameters. We separate relative errors for the image and motion. Plots of these errors against iteration can be seen in Fig.~\ref{fig:2D_SuperRes_RE} for the problem with $2$\% added noise. The corresponding resulting image for LAP for the $2$\% error case can be seen in Fig.~\ref{fig:2D_SuperRes_Images}. Lastly, a table of relevant values including the average number of iterations, minimum errors for image and motion, matrix-vector multiplications, and CPU timings for the methods taken over $10$ different realizations of the problem for all three noise levels is in Table~\ref{tab:2D_SuperRes_Table}. 

For direct regularization using the discrete gradient operator, the solutions for all three methods are comparable in terms of the relative error for the motion, with LAP and BCD slightly outperforming VarPro for the relative error of the recovered images. This is a direct result of the element-wise bound constraints on the resultant images using these methods. Furthermore, these solutions are superior to those for all three methods using hybrid regularization or the identity operator, suggesting that this is a more appropriate regularizer for this problem. LAP with the discrete gradient operator recovers the most accurate reconstructed image of the three methods and achieves better or comparable recovery of the motion parameters. This is observable in the relative error plots for the $2$\% added noise case in Fig.~\ref{fig:2D_SuperRes_RE}, and for the problem with all three noise levels in Table~\ref{tab:2D_SuperRes_Table}. We can also see from the relative error plots that LAP tends to recover the correct motion parameters earlier in the Gauss--Newton iterations than either BCD or VarPro. In terms of method cost, both the LAP and BCD iterations cost significantly less in terms of time and matrix-vector multiplications than those of VarPro, resulting in faster CPU times and fewer matrix-vector multiplies for the entire optimization. However, while BCD is also relatively cheap in terms of matrix-vector multiplies and CPU time, LAP outperforms it in terms of solution quality. Overall, the LAP approach compares favorably to VarPro and BCD for this example in terms of both the resulting solutions and cost, outperforming both methods.

\begin{table}[t]\small
    \begin{center}
        \footnotesize
        \begin{tabular}{|c|c|c|c|c|c|c|}
            \hline
            \multirow{2}{*}{\rottext{\hspace{4mm} 1\% Noise}}
            & & Iter. & Rel.~Err.~$\bfx$ & Rel.~Err.~$\bfw$  & MatVecs. & Time(s) \\ \cline{2-7}
& LAP + \texttt{HyBR}     & {\bf 11.1}    & 5.34e-2         & 1.82e-2        & 93.4        & {\bf 14.3} \\    
& LAP + $\nabla_h$         & 13.1         & {\bf 3.76e-2}     & {\bf 1.78e-2}     & {\bf 62.2}     & 15.4 \\    
& VarPro + $\bfI$         & 25.7         & 8.60e-2         & 1.84e-2         & 554.0         & 60.1 \\
& VarPro + $\nabla_h$     & 21.1         & 4.12e-2         & 1.81e-2         & 462.0         & 56.0 \\
& BCD + \texttt{HyBR}     & 12.0         & 5.93e-2         & 2.03e-2         & 68.1          & 22.9 \\
& BCD + $\nabla_h$     & 29.7         & 3.95e-2         & 1.79e-2         & 89.6         & 55.4 \\ \hline     
            \hline
                    \multirow{2}{*}{\rottext{\hspace{4mm} 2\% Noise}}
            & & Iter. & Rel.~Err.~$\bfx$ & Rel.~Err.~$\bfw$  & MatVecs. & Time(s) \\ \cline{2-7}
& LAP + \texttt{HyBR}     & {\bf 12.0}    & 6.38e-2         & 1.86e-2         & 77.0         & {\bf 12.5} \\    
& LAP + $\nabla_h$         & 15.9         & {\bf4.38e-2}     & {\bf 1.81e-2}    & 66.4         & 15.9 \\
& VarPro + $\bfI$         & 29.6         & 1.03e-1         & 1.86e-2         & 632.0         & 64.7 \\    
& VarPro + $\nabla_h$     & 21.4         & 5.05e-2         & 1.85e-2         & 468.0         & 54.9 \\
& BCD + \texttt{HyBR}     & 13.4         & 7.10e-2         & 2.15e-2        & {\bf 53.1}     & 25.0 \\
& BCD + $\nabla_h$     & 29.9         & 4.55e-2         & 1.82e-2         & 91.7         & 57.0 \\ \hline 
            \hline
            \multirow{2}{*}{\rottext{\hspace{4mm} 3\% Noise}}
            &     & Iter. & Rel.~Err.~$\bfx$ & Rel.~Err.~$\bfw$  & MatVecs. & Time(s) \\ \cline{2-7}
& LAP + \texttt{HyBR}     &  {\bf 12.3}    & 7.54e-2      & 2.28e-2         & 68.6        & {\bf 16.3} \\    
& LAP + $\nabla_h$         &  18.5        & {\bf 5.34e-2}    & 2.22e-2      & 74.2         & 22.0 \\
& VarPro + $\bfI$         &  34.4        & 1.27e-1         & 2.32e-2         & 728.0         & 86.7 \\    
& VarPro + $\nabla_h$     &  23.9        & 6.17e-2        & {\bf 2.13e-2}    & 518.0        & 72.5 \\
& BCD + \texttt{HyBR}     &  15.0        & 7.86e-2        & 2.32e-2        & {\bf 52.9}    & 31.9 \\
& BCD + $\nabla_h$     &  29.7         & 5.42e-2         & 2.24e-2        & 83.2        & 63.9 \\ \hline 
         \end{tabular}
    \end{center}
    \caption{This table shows data for the 2D super-resolution for multiple values of added Gaussian noise. The columns from left to right give the stopping iteration, the relative error of the solution image, the relative error of the solution motion, number of matrix-vector multiplications during optimization, and time in seconds using \texttt{tic} and \texttt{toc} in \texttt{MATLAB}. All values are averages taken from $10$ instances with different motion parameters, initial guesses, and noise realizations.
    The best results for each column is bold-faced.}
    \label{tab:2D_SuperRes_Table}
\end{table}

% new section
\subsection{Three-Dimensional Super Resolution} % (fold)
\label{sub:three_dimensional_example}

The next problem is a larger three-dimensional super-resolution problem. Again, we use a 3D MRI dataset provided in \texttt{FAIR}~\cite{Modersitzki2009} to construct a super-resolution problem with a known ground truth image and motion parameters. The ground truth image (resolution $160 \times 96 \times 144$) is used to generate 128 frames of low-resolution test data (resolution $40 \times 24 \times 32$). Each frame of data is shifted and rotated by a random 3D rigid body transformation, after which it is downsampled using block averaging. Lastly, Gaussian white noise is added to each low-resolution data frame in the same way as for the two-dimensional super-resolution problem, and we run the problem for $\mu = 1$\%, $2$\%, and $3$\% added noise per data frame. The resulting optimization problem has $2,212,608$ unknowns, $2,211,840$ for the image and $768$ for the motion parameters. The data has dimension $5,898,240$. The formulation of the problem is identical to that of the two-dimensional super-resolution problem with appropriate corrections for the change in dimension. The imaging operator $\bfK$ in the three-dimensional example is block diagonal with $128$ identical down-sampling matrices $K$ which relate the high-resolution image to the low-resolution data by block averaging.

The initial guess for the three-dimensional problem is generated using the same strategy as in the two-dimensional case. To generate a guess for the motion parameters, we register all frames onto the first frame (thus solving $127$ rigid registration problems.) This gives an initial guess for the motion with approximately $3$\% relative error. Using this initial guess, we then solve a linear least-squares problem to obtain an initial guess for the image with a relative error of around $13$\%. We note that this is a poorer initial guess when compared to the true solution, than the one obtained for the two-dimensional super-resolution example. This may impact the quality of the solution obtainable for examples with large amounts of noise. 

For this example, we test LAP, VarPro, and BCD using only the discrete gradient regularizer due to its better performance in the 2D example. We use identical parameters to the 2D problem for the regularization parameter, bound constraints on the image variables, and accuracy of the iterative LSQR. The lone difference in the problem setup from the 2D case is in solving the linear system \eqref{eq:VarPro1} for VarPro.  Instead of the $20$ fixed LSQR iterations from the 2D case, we run a fixed number of $50$ iterations to achieve the required accuracy for the larger system. For all three methods, the reduced system in the motion is solved using Cholesky factorization on the normal equations.

\begin{figure}[t]
    \begin{center}
        \begin{tabular}{c}
            \includegraphics[width=.95\textwidth]{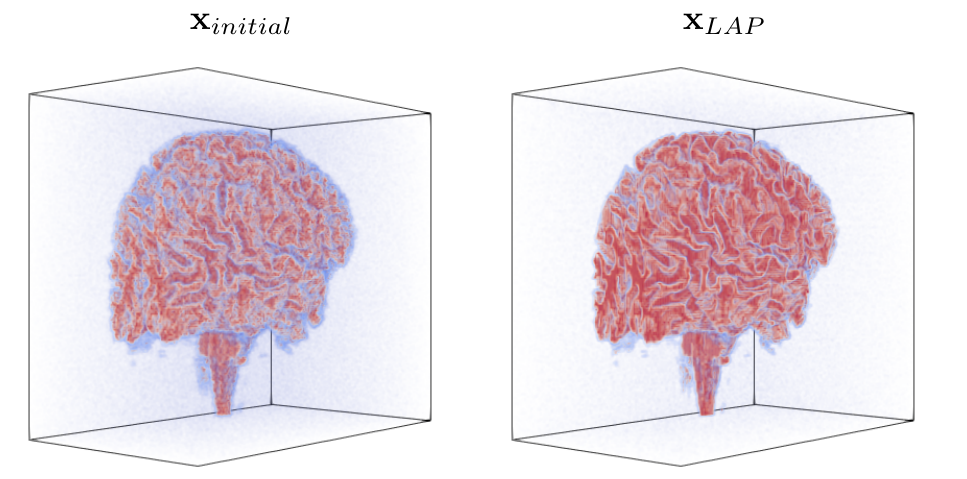} 
        \end{tabular}
    \end{center}
    \caption{This figure shows 3D volume renderings of the reconstructed images using LAP for the 3D super resolution problem with 2\% noise. The volume on the left shows the initial guess, $\bfx_{initial}$. On the right is the reconstructed solution using LAP, $\bfx_{LAP}$ with the discrete gradient regularizer.}
    \label{fig:3D_SuperRes_Images}
\end{figure}

\begin{figure}[t]
    \begin{center}
        \includegraphics[width=\textwidth]{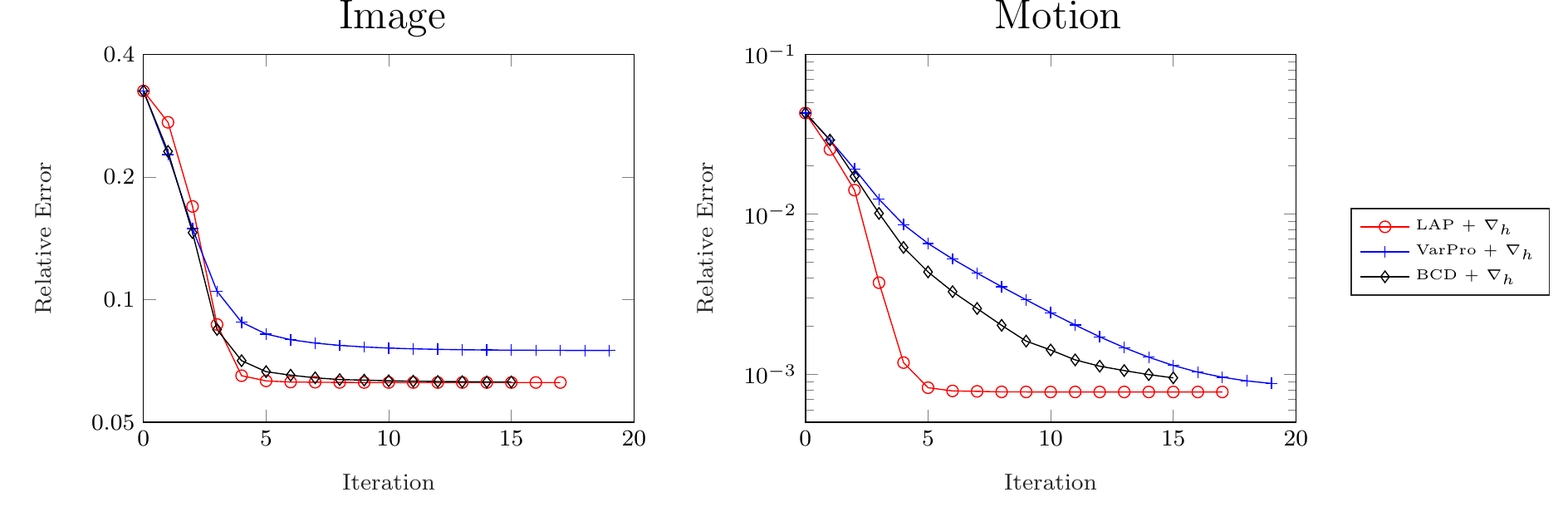}
    \end{center}
    \caption{This figure plots the relative errors for both the reconstructed image and the motion parameters for the 3D super resolution for $2$\% added noise. LAP succeeds in capturing the correct motion parameters in fewer iterations than VarPro and BCD and in recovering images of comparable quality.}
    \label{fig:3D_SuperRes_RE}
\end{figure}

The resulting solution image for LAP and the relative error plots for all three methods for the images and motion parameters can be found in Figs.~\ref{fig:3D_SuperRes_Images} and~\ref{fig:3D_SuperRes_RE}, respectively. Like the 2D super-resolution example, LAP converges faster to the motion parameters in the early iterations than BCD and VarPro and succeeds in reaching lower relative errors for both the recovered image and motion parameters. Again, VarPro's iterations are far more expensive in terms of matrix-vector multiplications and CPU time as seen in Table~\ref{tab:3D_SuperRes_Table}. BCD performs similarly with LAP in terms of the reconstructed image, but it does less well at recovering the motion parameters and is slightly more expensive in terms of CPU time due to a higher number of function calls.

\begin{table}[t]\small
    \begin{center}
        \footnotesize
        \begin{tabular}{|c|c|c|c|c|c|c|}
            \hline
            \multirow{2}{*}{\rottextt{ 1\% Noise}}
             & & Iter. & Rel.~Err.~$\bfx$ & Rel.~Err.~$\bfw$  & MatVecs. & Time(s) \\ \cline{2-7}
& LAP + $\nabla_h$         & 14.1         & {\bf 6.23e-2}         & {\bf 8.65e-4}         & {\bf 22.1}     & {\bf 2.30e3} \\    
& VarPro + $\nabla_h$     & 16.4         & 7.21e-2             & 9.55e-4             & 1.79e3         & 1.23e4 \\
& BCD + $\nabla_h$     & {\bf 11.3}     & 6.25e-2             & 1.40e-3             & 25.7         & 4.26e3 \\ \hline     
            \hline
            \multirow{2}{*}{\rottextt{ 2\% Noise}}
            & & Iter. & Rel.~Err.~$\bfx$ & Rel.~Err.~$\bfw$  & MatVecs. & Time(s) \\ \cline{2-7}
& LAP + $\nabla_h$         & 15.3         & {\bf 6.35e-2}     & {\bf 9.27e-4}     & {\bf 21.0}     & {\bf 2.15e3} \\    
& VarPro + $\nabla_h$     & 18.0         & 7.57e-2         & 1.00e-3         & 1.95e3         & 1.22e4 \\
& BCD + $\nabla_h$     & {\bf 12.2}     & 6.37e-2         & 1.59e-3         & 26.2         & 3.75e3 \\ \hline     
            \hline
            \multirow{2}{*}{\rottextt{ 3\% Noise}}
            & & Iter. & Rel.~Err.~$\bfx$ & Rel.~Err.~$\bfw$  & MatVecs. & Time(s) \\ \cline{2-7}
& LAP + $\nabla_h$         & 16.2         & {\bf 6.48e-2}     & {\bf 8.97e-4}     & {\bf 20.3}     & {\bf 2.94e3} \\    
& VarPro + $\nabla_h$     & 17.2         & 8.08e-2         & 9.52e-4         & 1.88e3         & 1.07e4 \\
& BCD + $\nabla_h$     & {\bf 11.7}     & 6.50e-2         & 1.53e-3         & 25.2         & 4.21e3 \\ \hline     
         \end{tabular}
    \end{center}
    \caption{This table presents data for the solution of the 3D super resolution for multiple values of added Gaussian noise. The columns from left to right give the stopping iteration, relative error of the solution image, relative error of the solution motion, number of matrix-vector multiplies during optimization, and time in seconds using \texttt{tic} and \texttt{toc} in \texttt{MATLAB}. All values are averages taken from $10$ separate problems with different motion parameters, initial guesses, and noise realizations.}
    \label{tab:3D_SuperRes_Table}
\end{table}

% new section
\subsection{MRI Motion Correction} % (fold)
\label{sub:MRI_motion_example}

The final test problem is a two-dimensional MRI motion correction problem.
The goal in this MRI application is to reconstruct a complex-valued MRI image from its Fourier coefficients that are acquired block-wise in a sequence of measurements.
Since the measurement process typically requires several seconds or minutes, in some cases the object being imaged moves substantially. Motion renders the Fourier samples inconsistent and --- without correction~--- results in artifacts and blurring in the reconstructed MRI image. To correct for this, one can instead view the collected MRI data as a set of distinct, complex-valued Fourier samplings, each measuring some portion of the Fourier domain and subject to some unknown motion parameters. The resulting problem of recovering the unknown motion parameters for each Fourier sampling and combining them to obtain a single motion-corrected MRI image fits into the coupled imaging framework presented in this paper. 

The forward model for this problem was presented by Batchelor {\em et al.}~\cite{BatchelorEtAl2005}. In their formulation, our imaging operator $\bfK$ is again block diagonal with diagonal blocks $K_k$ for $k = 1,2,\ldots,N$  given by
\begin{equation*}
K_k = \bfA_k \mathcal{F} \bfS.
\end{equation*}
\noindent Here, $\bfS$ is a complex-valued block rectangular matrix containing the given coil sensitivities of the MRI machine, $\mathcal{F}$ is block diagonal with each block a two-dimensional Fourier transform (2D FFT), and $\bfA_k$ is a block diagonal matrix with rectangular blocks containing selected rows of the identity corresponding to the Fourier sampling for the $k$th data observation. As with the other examples, the imaging operator $\bfK$ is multiplied on the right by the block rectangular matrix $\bfT$ with $T(\bfy(\bfw_k))$ blocks modeling the motion parameters of each Fourier sampling. We note that the cost of matrix-vector multiplications by this imaging operator is dominated by the 2D FFTs in block $\mathcal{F}$, and we note that for $32$ receiver coils, the cost of a single matrix-vector multiplication will require $32N$ 2D FFTs of size $128 \times 128$ where $N$ is the number of Fourier samplings in the data set. Additionally, we note that the presence of these FFT matrices prevents us from explicitly storing the matrix and necessitates passing it as a function call for all of the methods. This also applies to the Jacobian with respect to the image, $\bfJ_x$. However, because it is relatively small in size, $\bfJ_w$ can still be computed and stored explicitly.

We use the dataset provided in the \texttt{alignedSENSE} package~\cite{CorderoEtAl2016} to set up an MRI motion correction with a known true image and known motion parameters. To this end, we generate noisy data by using the forward problem~\eqref{eq:fwd}. The ground truth image with resolution $128 \times 128$ is rotated and shifted by a random 2D rigid body transformation. The motion affected image is then observed on $32$ sensors with known sensitivities. Each of these $32$ observations is then sampled in Fourier space. For our problem, each sampling corresponds to $1/16$ of the Fourier domain, meaning that $N = 16$ samplings (each with unknown motion parameters) are needed to have a full sampling of the whole space. We sample using a Cartesian parallel two-dimensional sampling pattern~\cite{CorderoEtAl2016}. Noise is added to the data using the formula
\begin{equation*}
\begin{array}{lll}
\bfd = \bar{\bfd} + \epsilon & \text{ where } & \epsilon = \mu \frac{\|\bar{\bfd} \|_\infty}{\| \bfn \|_2} \bfn.
\end{array} 
\end{equation*}
\noindent Here, $\bar{\bfd}$ is the noise free data, $\mu$ is the percentage of noise added, and $\bfn$ is a complex-valued vector where the entries of $\text{Re}(\bfn)$ and $\text{Im}(\bfn)$ are normally-distributed random numbers with mean $0$ and standard deviation $1$. We run the problem for $\mu = 5$\%, $10$\%, and $15$\% added noise. The resulting data has dimension $\textstyle{\frac{ 128 \times 128 \times 32}{16}} \times 16 = 524,288$. The MRI motion correction optimization problem fits within the framework in Eq.~\ref{eq:optProb} and has $16,432$ unknowns corresponding to $16,384$ for the image and $48$ for the motion.

As with previous examples, we solve the MRI motion correction problem using LAP, VarPro, and BCD. The setup and parameters are similar to those of the 2D super-resolution. We use the discrete gradient and \texttt{HyBR} as regularization options for LAP and BCD, and for VarPro we regularize using the discrete gradient operator and the identity operator. For the non-hybrid regularizers, we fix $\alpha = 0.01$. For the least-squares problems in the image variables for LAP and BCD, we solve using LSQR with an identical tolerance to the super-resolution examples. For VarPro, we use LSQR with a tolerance of $\num{e-8}$ or a maximum of $100$ iterations to solve \eqref{eq:VarPro1} to maintain accuracy in the gradient. Cholesky factorization on the normal equations is used for the lower-dimensional solves in the motion. No bound constraints are applied to any of the methods for this example because element-wise bound constraints on the real or imaginary parts of the complex-valued image variables will affect the angle (or phase) of the solution, which is undesirable.

\begin{figure}[t]
    \begin{center}
        \begin{tabular}{c}
            \includegraphics[width=.95\textwidth]{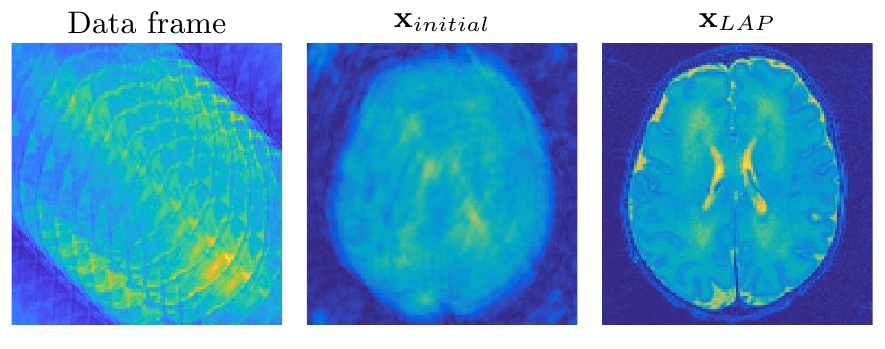} 
        \end{tabular}
    \end{center}
    \caption{This figure shows a single inverted data sampling (\textit{left}), the initial guess for the image (\textit{center}), and reconstructed image (\textit{right}) for the MRI motion correction problem with 10\% noise. The solution image shown is for LAP using the \texttt{HyBR} regularizer. Note that the MRI images are the modulus of the complex-valued images recovered.}
    \label{fig:MoCoMRI_Images}
\end{figure}

For an initial guess for the motion parameters, we start with $\bfw = 0$  (corresponding to a relative error of $100$\%). Using this initialization, we solve a linear least-squares problem to get an initial guess for the image. For $10$\% added Gaussian noise in the data, the initial guess for the image has a relative error of around $35$\%. We show the initial guess in Fig.~\ref{fig:MoCoMRI_Images}. 

\begin{figure}[t]
    \begin{center}
        \includegraphics[width=\textwidth]{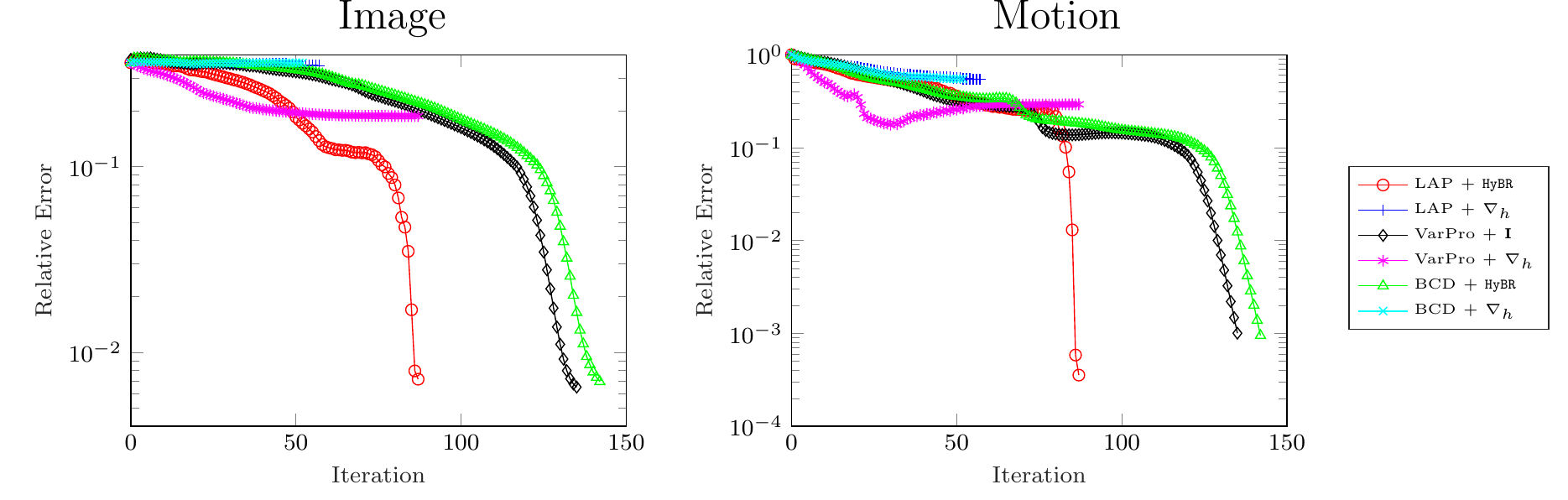}
    \end{center}
    \caption{This figure shows the relative errors for both the reconstructed image and the motion parameters for the MRI motion correction problem for three levels of noise. LAP with hybrid regularization achieves better reconstructions of the image and motion parameters for fewer iterations than either VarPro and BCD.}
    \label{fig:MoCoMRI_Plots}
\end{figure}

LAP, VarPro, and BCD provide fairly accurate reconstructions of both the image and motion parameters for quite large values of noise using either HyBR or the identity as a regularizer; see Figs.~\ref{fig:MoCoMRI_Images} and~\ref{fig:MoCoMRI_Plots}. This is likely due to the fact that the problem is not severely ill-posed and is highly over-determined (32 sensor readings for each point in Fourier space.) For this example, the hybrid regularization approach for LAP and BCD produces the best results, with LAP  requiring considerably fewer iterations. We remark that the best regularization from this problem differs from the super-resolution problems, which shows the importance of the flexibility that LAP offers for regularizing the image. The comparative speed of LAP is observable for the relative error plots for the problem with $10$\% noise and further evidenced in Table~\ref{MoCoMRI_Table} for all noise levels over $10$ separate realizations of the problem. For the gradient-based regularizer, all three methods do not recover the motion parameters accurately. We also note that the number of iterations and their cost is an important consideration for this problem. Because of the distance of the initial guess from the solution, this problem requires more iterations than the super-resolution examples. Additionally, the high number of 2D FFTs required for a single matrix-vector multiplication makes multiplications by the Jacobian $\bfJ_x$ and $\bfJ_w$ expensive. Table~\ref{MoCoMRI_Table} shows that LAP outperforms VarPro and BCD for both choices of regularizer by requiring fewer, cheaper iterations in terms of both time and matrix-vector multiplications. The difference in cost is most dramatic when compared with VarPro again due to the large number of FFTs required for a single matrix-vector multiplication and the large number of such multiplications required within each VarPro function call. For BCD and LAP, the number of matrix-vector multiplications is similar, but BCD requires more iterations for convergence. Overall, we see that LAP is a better choice for this problem and that it provides better reconstructions of both the image and motion in fewer, cheaper iterations.

\begin{table}[t]\small
    \begin{center}
        \footnotesize
        \begin{tabular}{|c|c|c|c|c|c|c|}
            \hline
             \multirow{2}{*}{\rottext{ \hspace{4mm} 5\% Noise}}
              &  & Iter. & Rel.~Err.~$\bfx$ & Rel.~Err.~$\bfw$  & MatVecs. & Time(s) \\ \cline{2-7}
                  & LAP + \texttt{HyBR}         & 76.6         & {\bf 3.55e-3}     & {\bf 1.59e-4}     & {\bf 3.18e2}     & 5.86e2 \\ 
                  & LAP + $\nabla_h$         & {\bf 76.2}     & 3.74e-3         & 1.67e-4         & 3.62e2         & {\bf 4.32e2}\\    
              & VarPro + $\bfI$         & 116.0         & 7.93e-2         & 4.32e-2         & 2.19e4         & 1.20e4 \\
             & VarPro + $\nabla_h$     & 115.9         & 7.92e-2         & 4.33e-2         & 2.19e4         & 1.52e4 \\
                 & BCD + \texttt{HyBR}         & 128.6         & 4.35e-2         & 2.27e-2         & 4.03e2         & 1.01e3 \\
             & BCD + $\nabla_h$         & 116.6         & 7.83e-2         & 4.20e-2         & 4.27e2         & 8.59e2 \\ \hline     
             
             \hline
             \multirow{2}{*}{\rottext{ \hspace{4mm}10\% Noise}}
              &  & Iter. & Rel.~Err.~$\bfx$ & Rel.~Err.~$\bfw$  & MatVecs. & Time(s) \\ \cline{2-7}
                  & LAP + \texttt{HyBR}         & 76.8         & {\bf 6.56e-3}     & {\bf 3.25e-4}    & {\bf 3.11e2}     & 6.52e2 \\ 
                  & LAP + $\nabla_h$         & {\bf 76.1}     & 6.88e-3         & 4.15e-4         & 3.55e2         & {\bf 5.19e2}\\    
              & VarPro + $\bfI$         & 116.0         & 8.08e-2         & 4.33e-2         & 2.19e4         & 1.12e4 \\
             & VarPro + $\nabla_h$     & 116.0         & 8.08e-2         & 4.32e-2         & 2.19e4         & 1.27e4 \\
                 & BCD + \texttt{HyBR}         & 128.2         & 4.53e-2         & 2.21e-2         & 3.89e2         & 1.24e3 \\
             & BCD + $\nabla_h$         & 116.0         & 8.02e-2         & 4.21e-2         & 4.12e2         & 1.05e3 \\ \hline 
             
             \hline
             \multirow{2}{*}{\rottext{ \hspace{4mm}15\% Noise}}
              &  & Iter. & Rel.~Err.~$\bfx$ & Rel.~Err.~$\bfw$  & MatVecs. & Time(s) \\ \cline{2-7}
                  & LAP + \texttt{HyBR}         & 77.3         & {\bf 9.49e-3} & {\bf 4.69e-2} & {\bf 3.07e2} & 5.02e2 \\ 
                  & LAP + $\nabla_h$         & {\bf 75.0}     & 1.61e-2         & 2.92e-2         & 3.40e2         & {\bf 3.61e2} \\    
              & VarPro + $\bfI$         & 115.8         & 8.29e-2         & 4.36e-2         & 2.18e4         & 9.98e3 \\
             & VarPro + $\nabla_h$     & 115.7         & 8.28e-2         & 4.35e-2         & 2.18e4         & 1.27e4 \\
                 & BCD + \texttt{HyBR}         & 127.0         & 4.80e-2         & 2.21e-2         & 3.47e2         & 1.06e3 \\
             & BCD + $\nabla_h$         & 129.0         & 8.20e-2         & 4.17e-2         & 3.99e2         & 8.17e2 \\ \hline 
             
         \end{tabular}
    \end{center}
    \caption{This table shows the results of LAP, VarPro, and BCD for solving the MRI motion correction example for multiple regularizers and varying levels of added noise. Averaged over $10$ realizations of the problem, the columns are stopping iteration, the relative error of the solution image, the relative error of the solution motion, number of matrix-vector multiplies during optimization, and time in seconds using \texttt{tic} and \texttt{toc} in \texttt{MATLAB}.  LAP outperforms the other methods in terms of solution quality, computational cost, and CPU time.}
    \label{MoCoMRI_Table}
\end{table}

% new section
\section{Summary and Conclusion} % (fold)
\label{sec:summary}

We introduce a new method, called Linearize And Project (LAP), for solving large-scale inverse problems with coupled blocks of variables in a projected Gauss--Newton framework. Problems with these characteristics arise frequently in applications, and we exemplify and motivate LAP using joint reconstruction problems in imaging that aim at estimating image and motion parameters from a number of noisy, indirect, and motion-affected measurements. By design, LAP is most attractive when the optimization problem with respect to one block of variables is comparably easy to solve. LAP is very flexible in the sense that it supports different regularization strategies, simplifies imposing equality and inequality constraints on both blocks of variables, and does not require (as in the case of VarPro) that the forward problem depends linearly on one set of variables. In our numerical experiments using four separable nonlinear least-squares problems, we showed that LAP is competitive and often superior to VarPro and BCD with respect to accuracy and efficiency. 

LAP is as general as alternating minimization methods such as Block Coordinate Descent. However while BCD ignores the coupling between the variable blocks when computing updates, LAP takes it into consideration. Thus, in our experiments LAP requires a considerably smaller number of iterations, matrix-vector multiplications, and CPU time than BCD to achieve  similar accuracy. 

LAP is not limited to separable nonlinear least-squares problems and thus more broadly applicable than VarPro. Since LAP projects after linearization, it provides the opportunity to freely choose which block of variables gets eliminated. For example, in our numerical examples in Sec.~\ref{sub:two_dimensional_example}--\ref{sub:MRI_motion_example}, LAP eliminates the parameters associated with the motion (which are of comparably small dimension) when computing the search direction in the projected Gauss--Newton scheme. Due to the robustness of Gauss--Newton methods, it suffices to solve the imaging problem iteratively to low accuracy. By contrast, VarPro eliminates the image variables that enter the residual in a linear way. While this leads to a small-dimensional nonlinear optimization problem for the motion, each iteration requires solving the imaging problem to relatively high accuracy to obtain reliable gradient information; see also Sec~\ref{sub:varpro}. This can be problematic for large and ill-posed imaging problems, and thus 
LAP can in some cases reduce runtimes by an order of magnitude; see Tables~\ref{tab:3D_SuperRes_Table} and \ref{MoCoMRI_Table}.

It is worth noting that the key step of LAP, which is the projection onto one variable block when solving the approximated Newton system, is a block elimination, and the reduced system corresponds to the Schur-complement. 

To allow for comparison with VarPro we focussed on separable nonlinear least-squares problems. In future work, we will study the performance of LAP to solve general coupled nonlinear optimization problems.

% new section
\section{ACKNOWLEDGEMENT} % (fold)
\label{sec:acknowledgements}
This work is supported by Emory's University Research Committee and National Science Foundation (NSF) awards DMS 1522760 and DMS 1522599.

% new section
\bibliographystyle{abbrv}
\bibliography{2016-CoupledSuperRes.bib}
\end{document}